\documentclass[11pt]{article}
\usepackage{amsfonts,amscd,amssymb}
\newtheorem{definition}{Definition}[section]

\newtheorem{proposition}[definition]{Proposition}
\newtheorem{corollary}[definition]{Corollary}
\newtheorem{remark}[definition]{Remark}

\newtheorem{theorem}[definition]{Theorem}
\newtheorem{example}[definition]{Example}

\newcommand{\nat}{\mbox{$\;\natural \;$}}

\def\rawo\lonra{\longrightarrow}

\def\ot{\otimes}

\newcommand{\eqref}[1]{(\ref{eq:#1})}

\begin{document}
\title{Some bialgebroids constructed by Kadison and Connes-Moscovici 
are isomorphic
\thanks{Research partially supported by the EC programme LIEGRITS, 
RTN 2003, 505078, and by the project ``New techniques in  
Hopf algebras and graded ring theory'' of the Flemish and Romanian 
Ministries of Research.}}
\author
{Florin Panaite\\
Institute of Mathematics of the 
Romanian Academy\\ 
PO-Box 1-764, RO-014700 Bucharest, Romania\\
e-mail: Florin.Panaite@imar.ro\\
\and 
Freddy Van Oystaeyen\\
Department of Mathematics and Computer Science\\
University of Antwerp, Middelheimlaan 1\\
B-2020 Antwerp, Belgium\\
e-mail: Francine.Schoeters@ua.ac.be}
\date{}
\maketitle

\begin{abstract}
We prove that a certain bialgebroid introduced recently by Kadison is 
isomorphic to a bialgebroid introduced earlier by Connes and Moscovici. 
At the level of total algebras, the isomorphism is a consequence of the 
general fact that an L-R-smash product over a Hopf algebra 
is isomorphic to a diagonal crossed product.  
\end{abstract}
%%%%%%%%%%%%%%%%%%%%%%%%%%%%%%%%%%%%%%%%%%%%%%%%%%%
\section{Introduction}
%%%%%%%%%%%%%%%%%%%%%%%%%%%%%%%%%%%%%%%%%%%%%%%%%%%%
Let $H$ be a Hopf algebra with bijective antipode, $A$ a left $H$-module 
algebra and denote as usual $A^e=A\ot A^{op}$. In \cite{kadison}, 
Kadison constructed a bialgebroid (in the sense of \cite{lu}) 
with base $A$ and having as total algebra a  
certain algebra structure on $A^e\ot H$, which we denote by 
$A^e\diamond H$. We first note that this algebra structure is actually  
an L-R-smash product $A^e\nat H$ as introduced in \cite{pvo}. An 
L-R-smash product is isomorphic to a diagonal crossed product as in 
\cite{hn1}, \cite{bpvo2}, hence $A^e\nat H\simeq A^e\bowtie H$. Finally, 
the diagonal crossed product  
$A^e\bowtie H$ is isomorphic to a certain algebra, which  
we denote by $A\odot H\odot A$, used by Connes and Moscovici in 
\cite{cm} (see also \cite{kr}) and which is also a bialgebroid over $A$. 
It turns out that the resulting algebra isomorphism is actually a  
bialgebroid isomorphism between $A^e\diamond H$ and $A\odot H\odot A$.\\
We establish also a universal property of $A^e\diamond H$ as a  
bialgebroid and we give an example (up to isomorphism) of the type 
$A\odot H\odot A$: the Cibils-Rosso algebra from \cite{cr}. 
%%%%%%%%%%%%%%%%%%%%%%%%%%%%%%%%%%%%%%%%%%%%%%%%%%%%%%%%%%%%%%
\section{The isomorphism}
%%%%%%%%%%%%%%%%%%%%%%%%%%%%%%%%%%%%%%%%%%%%%%%%%%%%%%%%%%%%%%%
\setcounter{equation}{0}
%%%%%%%%%%%%%%%%%%%%%%%%%%%%%%%%%%%%%%%%%%%%%%%%%%%%%%%%%%%%%%
We work over a ground field $k$. All algebras, linear spaces 
etc. will be over $k$; unadorned $\ot $ means $\ot_k$. Throughout $H$ will 
be a Hopf algebra with bijective antipode $S$. We use the following 
version of Sweedler's $\Sigma $-notation: $\Delta (h)=h_1\ot h_2$ 
for $h\in H$. \\[2mm]
Recall that the L-R-smash product over a cocommutative Hopf algebra 
was introduced in    
\cite{b1}, \cite{b2}, \cite{b3}, \cite{b4}, and generalized to an arbitrary 
Hopf algebra in \cite{pvo} as follows: if ${\cal A}$ is an  
$H$-bimodule algebra, 
the L-R-smash product ${\cal A}\nat H$ is the following 
algebra structure on ${\cal A}\ot H$: 
\begin{eqnarray*}
&&(\varphi \nat h)(\varphi '\nat h')=(\varphi \cdot h'_2)(h_1\cdot \varphi ')
\nat h_2h'_1, \;\;\;\forall \;\varphi , \varphi '\in {\cal A}, \;
h, h'\in H. 
\end{eqnarray*}
The diagonal crossed product ${\cal A}\bowtie H$ is the following algebra 
structure on ${\cal A}\ot H$, see \cite{hn1}, \cite{bpvo2}:
\begin{eqnarray*}
&&(\varphi \bowtie h)(\varphi '\bowtie h')=
\varphi (h_1\cdot \varphi '\cdot S^{-1}(h_3))
\bowtie h_2h', \;\;\;\forall \;\varphi , \varphi '\in {\cal A}, \;
h, h'\in H. 
\end{eqnarray*} 
\begin{proposition} (\cite{pvo}) ${\cal A}\nat H\simeq {\cal A}\bowtie H$ 
as algebras, the isomorphism being given as follows:
\begin{eqnarray*}
&&\nu :{\cal A}\bowtie H\rightarrow {\cal A}\nat H, \;\;\;\;
\nu (\varphi \bowtie h)=\varphi \cdot h_2\nat h_1, \\
&&\nu ^{-1}:{\cal A}\nat H\rightarrow {\cal A}\bowtie H, \;\;\;\;
\nu ^{-1}(\varphi \nat h)=\varphi \cdot S^{-1}(h_2)\bowtie h_1.
\end{eqnarray*}
\end{proposition}
Let now $A$ be a left $H$-module algebra, denote by $A^{op}$ the opposite 
algebra and by $A^e=A\ot A^{op}$. The bialgebroid introduced by Kadison in 
\cite{kadison} has as total algebra the following algebra structure on 
$A^e\ot H$, which we denote by $A^e\diamond H$:
\begin{eqnarray*}
&&(a\ot b\ot h)(a'\ot b'\ot h')=a(h_1\cdot a')\ot b'(S(h'_2)\cdot b)\ot 
h_2h'_1, 
\end{eqnarray*}
for all $a, a', b, b'\in A$ and $h, h'\in H$, where the 
multiplication on the second component is made in $A$ ({\it not} in 
$A^{op}$).  
On the other hand, the opposite algebra $A^{op}$ becomes a right 
$H$-module algebra with action $a\cdot h=S(h)\cdot a$, hence 
$A^e=A\ot A^{op}$ is an $H$-bimodule algebra with actions 
$h\cdot (a\ot b)\cdot h'=h\cdot a\ot b\cdot h'$,  
for all $a, b\in A$ and $h, h'\in H$, so we may consider the L-R-smash  
product $A^e\nat H$. 
\begin{proposition}
$A^e\diamond H=A^e\nat H$.
\end{proposition}  
The bialgebroid introduced by Connes and Moscovici in \cite{cm} and 
further studied in \cite{kr} has as total algebra the following 
algebra structure on $A\ot H\ot A$, which we denote by $A\odot H\odot A$: 
\begin{eqnarray*}
&&(a\ot h\ot b)(a'\ot h'\ot b')=a(h_1\cdot a')\ot h_2h'\ot  
(h_3\cdot b')b, 
\end{eqnarray*}
for all $a, a', b, b'\in A$ and $h, h'\in H$. We may consider again the  
$H$-bimodule algebra $A^e$ as above and the diagonal crossed product 
$A^e\bowtie H$.  
\begin{proposition} \label{dia}
We have an algebra isomorphism     
\begin{eqnarray*}
&&A\odot H\odot A\simeq (A\ot A^{op})\bowtie H, \;\;a\ot h\ot b\mapsto   
(a\ot b)\bowtie h.
\end{eqnarray*}
\end{proposition}
\begin{corollary}\label{cucu}
$A^e\diamond H\simeq A\odot H\odot A$ as algebras, an explicit isomorphism 
being given by 
\begin{eqnarray*}
&&A^e\diamond H\rightarrow A\odot H\odot A, \;\;\;\;a\ot b\ot h\mapsto 
a\ot h_1\ot h_2\cdot b, \\
&&A\odot H\odot A\rightarrow A^e\diamond H, \;\;\;\;a\ot h\ot b\mapsto 
a\ot S(h_2)\cdot b\ot h_1. 
\end{eqnarray*} 
\end{corollary}
We recall from \cite{kadison}, \cite{cm}, \cite{kr} the structure of 
$A^e\diamond H$ and $A\odot H\odot A$ as bialgebroids over $A$:\\[2mm]
(1) source maps: $A\rightarrow A^e\diamond H$, $a\mapsto (a\ot 1)\ot 1$ and 
$A\rightarrow A\odot H\odot A$, $a\mapsto a\ot 1\ot 1$; \\
(2) target maps:  
$A\rightarrow A^e\diamond H$, $a\mapsto (1\ot a)\ot 1$ and  
$A\rightarrow A\odot H\odot A$, $a\mapsto 1\ot 1\ot a$; \\
(3) comultiplications: 
\begin{eqnarray*}
&&A^e\diamond H\rightarrow (A^e\diamond H)\ot _A(A^e\diamond H), \;\;\;
(a\ot b)\ot h\mapsto ((a\ot 1)\ot h_1)\ot _A ((1\ot b)\ot h_2), \\
&&A\odot H\odot A\rightarrow (A\odot H\odot A)\otimes _A (A\odot H\odot A), 
\;\;\;a\ot h\ot b\mapsto (a\ot h_1\ot 1)\ot _A(1\ot h_2\ot b);
\end{eqnarray*}
(4) counits: 
\begin{eqnarray*}
&&\varepsilon :A^e\diamond H\rightarrow A, \;\;\;\varepsilon ((a\ot b)\ot h)=
a(h\cdot b), \\
&&\varepsilon :A\odot H\odot A\rightarrow A, \;\;\;\varepsilon (a\ot h\ot b)=
a\varepsilon (h)b.
\end{eqnarray*}
\begin{theorem}\label{main}
$A^e\diamond H$ and $A\odot H\odot A$ are isomorphic as bialgebroids 
over $A$, where the isomorphism is given as in Corollary \ref{cucu} at the 
level of total algebras and as identity at the level of base algebras. 
\end{theorem}
\begin{remark}{\em 
Assume that $H$ is involutive, i.e. $S^2=id$. In this case $A^e\diamond H$ 
and $A\odot H\odot A$ have antipodes making them Hopf algebroids 
(in the sense of \cite{bohm}), namely (see \cite{kadison}, \cite{kr}):
\begin{eqnarray*}
&&A^e\diamond H\rightarrow A^e\diamond H, \;\;\;(a\ot b)\ot h\mapsto  
(b\ot a)\ot S(h), \\
&&A\odot H\odot A\rightarrow A\odot H\odot A, \;\;\;
a\ot h\ot b\mapsto S(h_3)\cdot b\ot S(h_2)\ot S(h_1)\cdot a.
\end{eqnarray*}
Then one can see that the isomorphism from Theorem \ref{main} commutes 
with these antipodes, i.e. it is a {\it strict} isomorphism  of 
Hopf algebroids (in the terminology of \cite{bohm}).}
\end{remark} 
\begin{example}{\em 
Assume that $H$ is finite dimensional. Obviously $H^*$, the linear dual 
of $H$, is a left $H\ot H^{op}$-module algebra, with action 
$(h\ot h')\cdot f=h\rightharpoonup f\leftharpoonup h'$, for all 
$h, h'\in H$ and $f\in H^*$, where $\rightharpoonup $ and 
$\leftharpoonup $ are the regular actions of $H$ on $H^*$ given by 
$(h\rightharpoonup f)(h')=f(h'h)$ and $(f\leftharpoonup h')(h)=f(h'h)$. 
Hence, we can consider the bialgebroid $H^*\odot (H\ot H^{op})\odot H^*$ 
over $H^*$; by Proposition \ref{dia} we have 
$H^*\odot (H\ot H^{op})\odot H^*\simeq (H^*\otimes H^{* op})\bowtie 
(H\ot H^{op})$ as algebras. The algebra $Z=(H^*\ot H^{* op})\bowtie 
(H\ot H^{op})$ was considered before in \cite{pan}, it is isomorphic to the 
Cibils-Rosso algebra $X$ from \cite{cr} having the property that modules 
over it coincide with $H^*$-Hopf bimodules. Consequently, we obtain that 
the Cibils-Rosso algebra is a bialgebroid over $H^*$.}
\end{example}
%%%%%%%%%%%%%%%%%%%%%%%%%%%%%%%%%%%%%%%%%%%%%%%%%%%%%%%%%%%%%%%%
\section{A universal property of $A^e\diamond H$ as bialgebroid}
%%%%%%%%%%%%%%%%%%%%%%%%%%%%%%%%%%%%%%%%%%%%%%%%%%%%%%%%%%%%%
\setcounter{equation}{0}
%%%%%%%%%%%%%%%%%%%%%%%%%%%%%%%%%%%%%%%%%%%%%%%%%%%%%%%%%%%%%%
Kadison found in \cite{kadison} a universal property of $A^e\diamond H$ as 
algebra. We first establish an equivalent formulation, emphasizing the 
presence of the algebra $A\ot A^{op}$: 
\begin{proposition} \label{uni1}
Denote by $i:A\ot A^{op}\rightarrow (A\ot A^{op})\diamond H$, 
$i(a\ot b)=(a\ot b)\ot 1$ and $j:H\rightarrow (A\ot A^{op})\diamond H$, 
$j(h)=(1\ot 1)\ot h$, the canonical algebra inclusions. If $R$ is an 
associative algebra and $u:A\ot A^{op}\rightarrow R$ and 
$v:H\rightarrow R$ are algebra maps such that 
\begin{eqnarray*}
&&v(h_1)u(a\ot S(h_2)\cdot b)=u(h_1\cdot a\ot b)v(h_2),
\end{eqnarray*}
for all $a, b\in A$ and $h\in H$, there exists a unique algebra map 
$\omega :(A\ot A^{op})\diamond H\rightarrow R$ such that 
$\omega \circ i=u$ and $\omega \circ j=v$. The map $\omega $ is given by 
$\omega ((a\ot b)\ot h)=u(a\ot h_2\cdot b)v(h_1)$.
\end{proposition} 
Recall from \cite{lu} that $A\ot A^{op}$ is a bialgebroid over $A$, with 
source map $\alpha :A\rightarrow A\ot A^{op}$, $a\mapsto a\ot 1$, target 
map $\beta :A\rightarrow A\ot A^{op}$, $a\mapsto 1\ot a$, coproduct 
$\Delta :A\ot A^{op}\rightarrow (A\ot A^{op})\ot _A(A\ot A^{op})$, 
$\Delta (a\ot b)=(a\ot 1)\ot _A(1\ot b)$ and counit 
$\varepsilon :A\ot A^{op}\rightarrow A$, $\varepsilon (a\ot b)=ab$. With 
respect to this bialgebroid structure on $A\ot A^{op}$, we have a  
universal property of $(A\ot A^{op})\diamond H$ as bialgebroid:
\begin{theorem}
Let $R$ be as in Proposition \ref{uni1} and assume that $R$ is moreover a 
bialgebroid over $A$, $u$ is a morphism of bialgebroids (with identity  
map at the level of base algebras) and $v$ is a morphism of bialgebroids  
(with canonical inclusion $k\rightarrow A$ at the level of base algebras). 
Then $\omega $ is a morphism of bialgebroids (with identity map at the level 
of base algebras). 
\end{theorem}
The morphism of bialgebroids $(A\ot A^{op})\diamond H\rightarrow 
A\odot H\odot A$ from Theorem \ref{main} can be reobtained easily using this  
universal property.        
%%%%%%%%%%%%%%%%%%%%%%%%%%%%%%%%%%%%%%%%%%%%%%%%%%%%%%%%%%%%

\end{document}